\documentclass{article}

\newcommand{\R}{I\!\! R}
\newcommand{\C}{I\!\! C}
\newcommand{\Z}{I\!\! Z}
\newcommand{\N}{I\!\! N}

\begin{document}

Tsemo Aristide

The Abdus Salam Internation Centre for Theoretical Physics

Strada Costiera, 11

Trieste, Italy.

email: tsemo@ictp.trieste.it

\bigskip

\bigskip

{\centerline{\bf Cohomologie non ab\'elienne d'ordre sup\'erieur
et applications.}}

\bigskip

\bigskip

\bigskip

{\centerline{\bf Abstract.}}

{\it We propose in this paper a theory of higher non abelian
cohomology without using the notion of $n-$category. We use  this
theory to study compositions series of affine manifolds and
cohomology of manifolds.}

\bigskip

\bigskip

{\bf 0. Introduction.}

\bigskip

 Une vari\'et\'e affine $(M,\nabla_M)$, est une vari\'et\'e
 diff\'erentiable $M$, munie d'une connexion $\nabla_M$ dont les
 tenseurs de courbure et de torsion sont nuls. La vari\'et\'e
 affine de dimension $n$, $(M,\nabla_M)$ est compl\`ete si et
 seulement si elle est le quotient de ${\R}^n$ par un sous-groupe
 $\Gamma_M$
 de $Aff({\R}^n)$ qui agit proprement et librement sur ${\R}^n$.
 L. Auslander a conjectur\'e que le groupe fondamental d'une
 vari\'et\'e affine compacte et compl\`ete est polycyclique. Dans
 [T4], on a conjectur\'e qu'un rev\^etement fini d'une vari\'et\'e affine
 compacte et compl\`ete est le domaine de
 d\'efinition d'une application affine non triviale. Ceci conduit
 \`a se poser le probl\`eme de classifier les s\'eries de
 compositions $(M_n,\nabla_{M_n})\rightarrow...\rightarrow
 (M_1,\nabla_{M_1})$ o\`u chaque application
 $f_i:(M_{i+1},\nabla_{M_{i+1}})\rightarrow (M_i,\nabla_{M_i})$
 est une fibration affine de domaine de d\'efinition compact et
 complet; ce qui signifie que $f_i$ est
 surjective, et chaque vari\'et\'e affine $(M_i,\nabla_{M_i})$
 est compacte et compl\`ete.
 La classification des fibrations affines
  s'effectue avec la th\'eorie des
 gerbes (Voir [T6]). Il est naturel de penser que la classification des
 s\'eries de compositions devrait se faire avec la notion de
 $n-$gerbe. La d\'efinition d'une notion de $n-$cat\'egorie
 n'\'etant pas bien comprise, on \'effectue la classification des
 s\'eries
 de compositions par le biais de la notion de tour gerb\'ee.

N\'eanmoins, rappelons qu'il est effectu\'e dans [T6] une
classifications des s\'eries de compositions o\`u chaque
application $f_i$ est affinement localement triviale (ceci veut
dire que l'holonomie d'une fibre ne d\'epend pas de celle-ci) en
utilisant une partie de l'axiomatique n\'ecessaire pour d\'efinir
la notion de $n-$cat\'egorie.

En fait le probl\`eme g\'en\'eral auquel se ram\`ene la
classification des s\'eries de compositions de vari\'et\'es
affines est celui de l'\'elaboration d'une th\'eorie de
cohomologie non commutative. Dans la deuxi\`eme partie de ce
travail, on \'elabore une th\'eorie de cohomologie non ab\'elienne
en utilisant  la notion de tour gerb\'ee afin d'\'eviter justement
de parler de $n-$cat\'egories.

Comme le mentionne Brylinski [Bry], les classes caract\'eristiques
sont utilis\'ees par plusieurs auteurs pour cr\'eer divers objets
math\'ematiques. C'est ainsi que Witten [Wi] utilise des classes
caract\'eristiques d'ordre $2$ afin d'\'etudier le polynome de
Jones. Il apparait donc n\'ecessaire de donner une description
g\'eom\'etrique des classes caract\'eristiques. Dans cette
optique, il faut donner une interpr\'etation g\'eom\'etrique de la
cohomologie enti\`ere d'une vari\'et\'e diff\'erentiable $M$. Le
groupe de cohomologie $H^2(M,{\Z})$ s'interpr\`ete comme
l'ensemble des classes d'\'equivalences des fibr\'es complexes sur
$M$ d'apr\`es la th\'eorie de Kostant Weil. Dans [Bry], il est
donn\'e une interpr\'etation de $H^3(M,{\Z})$ en termes de classes
d'\'equivalences des faisceaux de groupoides de Dixmier-Douady.
Notre th\'eorie permet d'interpr\'eter les groupes
$H^{n+2}(M,{\Z})$ en termes de classes d'\'equivalences de
$n-$tours gerb\'ees.

\medskip

 Voici le plan de notre travail

0. Introduction.

1. Rappels sur la notion de gerbe.

2. Tour gerb\'ee.

3. Suites spectrales et tours gerb\'ees.

4.  Application des tours gerb\'ees \`a la g\'eom\'etrie affine.

5.  Interpr\'etation de la cohomologie enti\`ere d'une
vari\'et\'e.

 \bigskip
\bigskip

{\bf 1. Rappels sur la notion de gerbe.}

\medskip

On rappelle dans cette partie la notion de gerbe pr\'esent\'ee
dans [Gi].

\medskip

{\bf D\'efinitions 1.1.}

- Soit $C$ une cat\'egorie,  un crible de $C$ est une partie $R$
des objets $Ob(C)$ de $C$ telle que pour toute fl\`eche
$m:X\rightarrow Y$, si $Y\in R$, alors $X$ appartient \`a $R$.

- Pour toute foncteur $f:E'\rightarrow E$, on note $R^f=\{X\in
Obj(E'), f(X)\in R\}$; $R^f$ est un crible.

\medskip

{\bf D\'efinition 1.2.}

Une topologie sur $C$ est une application qui a tout objet $S$
associe une partie non vide $J(S)$ de l'ensemble des cribles de
la cat\'egorie $E_S$ au-dessus de $S$, telle que

pour toute fl\`eche $f:T\rightarrow S$ de $C$ et tout crible $R$
appartenant \`a $J(S)$, $R^f$ appartient \`a $J(T)$.

Pour tout objet $S$ de $C$,  tout \'el\'ement $R$ de $J(S)$, et
tout
 crible $R'$ de $E_S$, $R'$ appartient \`a $J(S)$ d\`es que pour
tout objet $f:T\rightarrow S$ de $R$, on a $R'^f$ appartient \`a
$J(T)$.

Les \'el\'ements de $J(S)$ sont appel\'es les raffinements de $S$.

\medskip

{\bf D\'efinitions 1.3}

- Un faisceau sur $C$ est un foncteur contravariant $F$ de $C$
dans la cat\'egorie des ensembles tel que pour tout crible $R$ de
$S$, l'application
$$
F(S)\longrightarrow {{lim}{(F\mid R)}}
$$
soit bijective, o\`u $F\mid R$ est le pr\'efaisceau sur $R$
d\'efini par $(F\mid R)(f)= F(T)$ pour toute fl\`eche
$f:T\rightarrow S$ de $R$.

\medskip

- Soit $F:M\rightarrow B$ un foncteur, pour tout objet $S$ de $B$,
on appelle cat\'egorie fibre $F_S$ de $F$ en $S$, la
sous-cat\'egorie de $M$ dont les fl\`eches $m$ v\'erifient
$F(m)=Id_S$. Ceci implique que les objets de $F_S$ se projettent
sur $S$ par $F$. Soient $X$ et $Y$ deux objets de $F_S$, on
d\'esigne par $Hom_S(X,Y)$, l'ensemble des $m:X\rightarrow Y$ qui
se projettent sur $Id_S$.

\medskip

{\bf D\'efinitions 1.4.}

- Soient $F:M\rightarrow B$ un foncteur, $m:x\rightarrow y$ une
fl\`eche de $M$ et $f:T\rightarrow S$ sa projection par $F$, on
dit que $m$ est $B-$cart\'esienne, ou que $m$ est une image
inverse de $y$ par $f$, ou encore que $x$ est une image inverse de
$y$ si pour tout objet $z$ de $M_T$, l'application
$$
Hom_T(z,x)\longrightarrow Hom_f(z,y)
$$
$$
n\mapsto mn
$$
est bijective, o\`u $Hom_f(z,y)=\{n\in Hom(z,y),F(n)=f\}.$

\medskip

- Une cat\'egorie $B$ fibr\'ee est un foncteur $F:M\rightarrow B$
tel que pour toute fl\`eche $f:T\rightarrow S$ de $B$ et tout
objet $y$ de $F_S$, il existe une image inverse de $y$ par $f$, et
telle que le compos\'e de deux morphismes $B-$cart\'esiens est
$B-$cart\'esien.

\medskip

- On dira que la cat\'egorie est fibr\'ee en groupoides, si pour
tout diagramme

$$
x{\buildrel{f}\over{\longrightarrow}}
z{\buildrel{g}\over{\longleftarrow}} y
$$
dans $M$, au-dessus de

$$
U{\buildrel{\phi}\over{\longrightarrow}}W
{\buildrel{\psi}\over{\longleftarrow}}V,
$$
pour toute fl\`eche $m:U\rightarrow V$ telle que $\phi=\psi m$,
alors il existe un unique $h:x\rightarrow y$ tel que $f=gh$ et
$F(h)=m$.

Ceci implique que l'image inverse d'une fl\`eche est unique \`a
isomorphisme pr\`es.

Soit $f:x\rightarrow y$ une fl\`eche de $M$ au-dessus de
$\phi:U\rightarrow V$ de $B$. On posera $x=\phi^*(y)$. Les
foncteurs $(\phi\psi)^*$ et $\psi^*\phi^*$ sont canoniquement
isomorphes.

\bigskip

{\bf D\'efinition 1.5.}

Une cat\'egorie fibr\'ee $F\rightarrow E$ est dite scind\'ee si
elle est munie d'une application (le scindage) qui \`a toute
fl\`eche $f:T\rightarrow S$ de $E$ et \`a tout objet $x$ de $F_S$
associe une image inverse $x_f:x^f\rightarrow x$ de sorte que pour
toute suite
$U{\buildrel{g}\over{\rightarrow}}T{\buildrel{f}\over{\rightarrow}}S$
et tout objet $y$ de $F_S$, on a $(x^f)^g=x^{fg}$.

\bigskip

Soit un diagramme de cat\'egories
$$
\matrix{F && G\cr \downarrow f && \downarrow g\cr A &
{\buildrel{a}\over{\longrightarrow}} & B}
$$
On note $Hom_a(F,G)$ la sous-cat\'egorie de $Hom(F,G)$ dont les
objets sont les $u:F\rightarrow G$ tels que $gu=af$, les fl\`eches
\'etant les morphismes $m:u\rightarrow u'$ tels que $gm$ soit le
morphisme identique du foncteur $af$.

On note $Cart_a(F,G)$ la sous-cat\'egorie pleine de $Hom_a(F,G)$
dont les objets sont les $a-$foncteurs cart\'esiens c'est \`a dire
les foncteurs qui transforment tout morphisme cart\'esien en un
morphisme cart\'esien.

\medskip

Soient $E$ une cat\'egorie munie d'une topologie, et $F\rightarrow
E$ une cat\'egorie fibr\'ee, pour tout raffinement $R$ de $C$, on
a un foncteur restriction $Cart_E(E_S,F)\rightarrow Cart_E(R,F)$.

\medskip

 {\bf D\'efinition 1.6.}

Soit $E$ une cat\'egorie munie d'une topologie, un champ  est une
cat\'egorie fibr\'ee $F\rightarrow E$ telle que le foncteur
restriction $Cart_E(E_S,F)\rightarrow Cart_E(R,S)$ est une
\'equivalence de cat\'egories.

\bigskip

{\bf Remarque.}

Supposons que les produits fibr\'es existent dans $E$ ainsi qu'une
famille couvrante $(S_i\rightarrow S)_{i\in I}$ pour la topologie
de $E$.

Se donner un champ en groupoides $F\rightarrow E$ est \'equivalent
\`a se donner une cat\'egorie fibr\'ee en groupoides $F\rightarrow
E$ telle que

(i) Recollement des fl\`eches.

Pour tout objet $U$ de $E$, et tout objets $x$ et $y$ de $F_U$ le
foncteur de $F_U$ dans la cat\'egorie des ensembles qui \`a une
fl\`eche $f:V\rightarrow U$ associe  $Hom_V(f^*(x),f^*(y))$ est un
faisceau.

(ii) Recollement des objets

Soit $x_i$ un objet de $F_{S_i}$, et $t_{ij}$ une fl\`eche entre
la restriction de $x_j$ \` a $S_i\times_S S_j$ et celle de $x_i$
\`a $S_i\times_S S_j$, alors si  sur $S_i\times_S S_j\times_S S_k$
on a $t_{ik}=t_{ij}t_{jk}$, il existe un objet $x$ de $F_S$ dont
la restriction \`a $F_{S_i}$ est $x_i$.

\medskip

 Si de plus il existe une famille couvrante $(S_i\rightarrow S)_{i\in I}$
 telle que $F_{S_i}$ soit non vide, et pour tout objets $x$ et $y$
 de $F_{S_i}$, $Hom_{S_i}(x,y)$ est non vide (connexit\'e locale),
 on dit alors que la
 cat\'egorie fibr\'ee est une gerbe.

  Si de plus  le pr\'efaisceau $L_F$ d\'efini sur $E$ par
 $$
 U\longrightarrow Hom_U(F_U)
 $$
 est un faisceau en groupes qui commute avec les isomorphismes
 et les restrictions, on dit alors que $L_F$ est le lien de
 la gerbe. Dans la suite on consid\`erera uniquement des gerbes
 li\'ees.

\medskip

Dans la suite, $E$ sera une cat\'egorie, admettant les produits
fibr\'es finis, munie d'une topologie d\'efinie par la famille
couvrante $(S_i\rightarrow S)_{i\in I}$.

\bigskip

{\bf Classification des gerbes li\'ees.}

\medskip

 Dans ce paragraphe, on rappelle la classification des gerbes
au-dessus d'une cat\'egorie $C$ de lien donn\'e $L$.

\medskip

{\bf D\'efinitions 1.7.}

- Une gerbe $F\rightarrow E$ est triviale, si elle admet une
section. Si la cat\'egorie $E$ a un objet final $e$, cela signifie
que $F_e$ n'est pas vide.

- Deux gerbes $F\rightarrow E$ et $F'\rightarrow E$ de lien $L$
sont \'equivalentes si et seulement si il existe un isomorphisme
de cat\'egories fibr\'ees entre $F$ et $F'$ qui commute avec
l'action de $L$.

On note $H^2(E,L)$ l'ensemble des classes d'\'equivalences des
gerbes.

\medskip

Soit $(x_i)_{i\in I}$ une famille d'objets tels que $x_i$
appartient \`a $F_{S_i}$. Il existe une fl\`eche $u_{ij}$ entre
les restrictions de $x_j$ et $x_i$ \`a $F_{S_i\times_S S_j}$.
Notons $u^k_{ij}$ la restriction de $u_{ij}$ entre les
restrictions  respectives de $x_j$ et $x_i$ \`a $F_{S_{ijk}}$,
o\`u $S_{ijk}=S_i\times_S S_j\times_S S_k$.

\bigskip

{\bf Th\'eor\`eme. [Gi].}

{\it la famille $c_{ijk}=u^i_{jk}(u^j_{ik})^{-1}u^k_{ij}\in
L(S_{ijk})$ est un $2-$cocycle de Cech qui classifie la gerbe
$F\rightarrow E$ si le lien est ab\'elien.}

\bigskip

{\bf 2. Tour gerb\'ee.}

\medskip

Le but de cette partie est de g\'en\'eraliser la notion de gerbe
\`a celle de tour gerb\'ee. Cette notion de tour gerb\'ee
permettra ult\'erieurement de d\'efinir des groupes de
cohomologies non ab\'eliens d'ordre sup\'erieurs.

\medskip

La notion de tour gerb\'ee est d\'efinie r\'ecursivement \`a
partir de la notion de gerbe. Soit $E_1$ une cat\'egorie munie
d'une topologie d\'efinie par une famille couvrante
$(S_i\rightarrow S)_{i\in I}$. On suppose que les produits
fibr\'es existent dans $E_1$.

\medskip

{\bf D\'efinition 2.1.}

On appelle $2-$tour gerb\'ee, la donn\'ee
  d'une gerbe li\'ee
$F_1\rightarrow E_1$ de lien $L_1$ telle que

A tout objet $X(S_1)$ de $F_{S_1}$ o\`u $S_1$ est un objet
quelconque de $E_1$, on associe une gerbe li\'ee
$F_2(S_1,X(S_1))\rightarrow E_2(S_1,X(S_1))$ de lien $L_2$ un
faisceau sur $E_1$.  Ceci veut dire qu'il existe un isomorphisme
entre le groupe des automorphismes de la fibre d'un objet de
$E_2(S_1,X(S_1))$ qui se projettent sur l'identit\'e et les
sections d\'efinies sur $S_1$, d'un faisceau de $E_1$. Cet
isomorphisme commute aux morphismes entre objets et aux
restrictions. La classe d'isomorphisme de
$F_2(S_1,X(S_1))\rightarrow E_2(S_1,X(S_1))$ ne d\'epend pas de
$S_1$ et de $X(S_1)$.

  La propri\'et\'e
de connexit\'e suivante est v\'erifi\'ee:

Soient $S_1$ un objet de $E_1$, $T_1$, $T_2$,  $T_3$ des objets de
$F_{S_1}$, il existe une application

$$
*: Hom(T_1,T_2)\longrightarrow Hom(E_2(S_1,T_1),E_2(S_1,T_2))
$$
$$
u\longrightarrow u^*
$$
Soit $X$ un objet de $E_2(S_1,T_1)$, $U$ un objet de
$F_2(S_1,T_1)_{X}$, pour tout objet $V$ de
$F_2(S_1,T_2)_{u^*(X)}$, il existe un morphisme $u^*(U,V)$ entre
les gerbes $F_2(S_1,T_1)\rightarrow E_2(S_1,T_1)$ et
$F_2(S_1,T_2)\rightarrow E_2(S_1,T_2)$ au-dessus de $u^*$ qui
transforme $U$ en $V$. On fera souvent l'abus de noter $u^*$ et
$u^*(U,V)$ par $u^*$.

Soient $u:T_1\rightarrow T_2$, et $v:T_2\rightarrow T_3$ des
\'el\'ements respectifs de $Hom(T_1,T_2)$ et $Hom(T_2,T_3)$,
$(vu)^*$ coincide  avec $v^*u^*$ sur les objets et on suppose que
$(vu)^*$ coincide avec $v^*u^*$ modulo un morphisme de
$F_2(S_1,T_3)\rightarrow E_2(S_1,T_3)$. On a $(Id)^*=Id$ sur les
objets (mais pas forc\'ement sur les fl\`eches).

Soient $\phi:S\rightarrow T$ une fl\`eche de $E_1$ et $y$ un objet
de ${F_1}_{T}$. Consid\'erons une image inverse $m:x\rightarrow y$
de $\phi$. On a un morphisme
$$
m^*: E_2(T,y)\longrightarrow E_2(S,x).
$$
Soient $X$ un objet de $E_2(T,y)$, $U$ un objet de $F_2(T,y)_{X}$
et $V$ un objet de $F_2(S,x)_{m^*(X)}$, il existe un morphisme
entre les gerbes $F_2(T,y)\rightarrow E_2(T,y)$ et
$F_2(S,x)\rightarrow E_2(S,x)$ au-dessus de $m^*$ qui transforme
$U$ en $V$.

Soient $\psi:R\rightarrow S$ une autre fl\`eche de $E_1$ et
$n:z\rightarrow x$ une image inverse de $x$, les restrictions
$(mn)^*$ et $n^*m^*$ coincident sur les objets et on suppose
qu'elles coincident modulo un morphisme de $F_2(R,z)\rightarrow
E_2(R,z)$.

\medskip

Il est commode de supposer que les cat\'egories $E_1$ et
$E(S_1,X(S_1))$ sont $U-$petites pour d\'efinir l'ensemble
$\bigcup_{S_1\in Ob(E_1)}E(S_1,X(S_1))$.

\medskip

{\bf D\'efinition 2.2.}

Une $n+1-$tour gerb\'ee de base $E_1$ consiste en la donn\'ee
d'une gerbe $F_1\rightarrow E_1$ de lien $L_1$,

pour tout objets $S_1$ de $E_1$ et $S_1'$ de ${F_1}_{S_1}$ d'une
gerbe $F_2(S_1,S_1')\rightarrow E_2(S_1,S_1')$ de lien $L_2$,

pour tout quadruple $(S_1,S_1',S_2,S_2')$ o\`u $S_1$ est un objet
de $E_1$ et $S_1'$ un objet de ${F_1}_{S_1}$, $S_2$ est un objet
de $E_2(S_1,S_1')$ et $S_2'$ un objet de ${F_2(S_1,S_1')}_{S_2}$
d'une gerbe $F_3(S_1,S_1',S_2,S_2')\rightarrow
E_3(S_1,S_1',S_2,S_2')$ de lien $L_3$.

Supposons d\'efinie la famille de gerbes
$F_{i}(S_1,S_1',...,S_{i-1},S_{i-1}')\rightarrow
E_{i}(S_1,S_1'...,S_{i-1},S_{i-1}')$. Pour tout objet $S_i$ de
$E_{i}(S_1,S_1'...,S_{i-1},S_{i-1}')$ et $S_i'$ de
$F_i(S_1,S_1',..,S_{i-1},S_{i-1}')_{S_i}$ on se donne une gerbe
$F_{i+1}(S_1,S_1'...,S_{i},S_i')\rightarrow
E_{i+1}(S_1,S_1'...,S_{i},S_i')$ de lien $L_{i+1}$. En d'autres
termes il existe un isomorphisme entre le groupe des
automorphismes (qui se projettent sur l'identit\'e) de la fibre
d'un objet de $E_{i+1}(S_1,..,S_i,S_i')$ et les sections
d\'efinies sur $S_1$ d'un faisceau $L_{i+1}$ de $E_1$. On suppose
en outre que cet isomorphisme commute aux morphismes entre objets
et aux restrictions.

On suppose que la classe d'isomorphisme de
$F_i(S_1,..,S_{i-1})\rightarrow E_i(S_1,..,S_{i-1})$ ne d\'epend
pas de $S_1$,...,$S_{i-1}$.

On notera souvent de mani\`ere abusive la $i-$gerbe
$F_i(S_1,S_1',..,S_{i-1},S'_{i-1})\rightarrow
E_i(S_1,S_1',..,S_{i-1},S_{i-1}')$ par $E_2(S'_{i-1},L_i)$.

 On d\'efinit ainsi par r\'ecurence la famille de gerbes
 $F_{n+1}(S_1,S_1',..,S_n,S_n')\rightarrow
 E_{n+1}(S_1,S_1',..,S_n,S_n')$.

On suppose que la propri\'et\'e suivante de connexit\'e est
v\'erifi\'ee:

Soient $S_{i+1}$ un objet de $E_{i+1}(S_1,S_1',...,S_i,S_i')$ et
$T_1$, $T_2$ et $T_3$ des objets de
$F_{i+1}(S_1,S_1',..,S_i,S_i')_{S_{i+1}}$, on a une application
$$
Hom(T_1,T_2)\longrightarrow
Hom(E_{i+2}(S_1,..,S_{i+1},T_1),E_{i+2}(S_1,..,S_{i+1},T_2))
$$
$$
u\longrightarrow u^*
$$

Consid\'erons maintenant un objet $X$ de
$E_{i+2}(S_1,..,S_{i+1},T_1)$, $U$ un objet de
$F_{i+2}(S_1,..,S_{i+1},T_1)_X$ et $V$ un objet de
$F_{i+2}(S_1,..,S_{i+1},T_2)_{u^*(X)}$, il existe un morphisme
$u^*(U,V)$ entre les gerbes
$F_{i+2}(S_1,..,S_{i+1},T_1)\rightarrow
E_{i+2}(S_1,..,S_{i+1},T_1)$ et
$F_{i+2}(S_1,..,S_{i+1},T_2)\rightarrow
E_{i+2}(S_1,..,S_{i+1},T_2)$ au-dessus de $u^*$ qui transforme $U$
en $V$.

\medskip

On fera souvent l'abus de noter $u^*$ et $u^*(U,V)$ par $u^*$.

 Soient $u:T_1\rightarrow T_2$, et $v:T_2\rightarrow T_3$
des \'el\'ements respectifs de $Hom(T_1,T_2)$ et $Hom(T_2,T_3)$,
le morphisme $(uv)^*$ coincide avec $u^*v^*$ sur les objets et on
suppose qu'ils coincident  modulo un morphisme de
$F_{i+2}(S_1,..,S_{i+1},T_3)\rightarrow
E_{i+2}(S_1,..,S_{i+1},T_3)$.

On a $(Id)^*=Id$ sur les objets (mais pas forc\'ement sur les
fl\`eches).

Soient $\phi:S\rightarrow T$ une fl\`eche de
$E_{i+1}(S_1,S_1',..,S_i,S_i')$ et $x$ un objet de
${F_{i+1}}(S_1,S_1',..,S_i,S_i')_{T}$. Consid\'erons une image
inverse $m:y\rightarrow x$ de $\phi$. Est associ\'e \`a $m$ la
restriction $m^*$, morphisme entre les cat\'egories
$E_{i+2}(S_1,..,T,x)$ et $E_{i+2}(S_1,..,S,y)$.

Consid\'erons maintenant un objet de $X$ de $E_{i+2}(S_1,..,T,x)$,
$U$ un objet de $F_{i+2}(S_1,..,T,x)_{X}$ et $V$ un objet de
$F_{i+2}(S_1,..,S,y)_{m^*(X)}$. Il existe un morphisme $m^*(U,V)$
entre les gerbes $F_{i+2}(S_1,..,T,x)\rightarrow
E_{i+2}(S_1,..,T,x)$ et $F_{i+2}(S_1,..,S,y)\rightarrow
E_{i+2}(S_1,..,S,y)$ au-dessus de $m^*$ qui transforme $U$ en $V$.

Soit $\psi:U\rightarrow S$ une autre fl\`eche de
$E_{i+1}(S_1,S_1',..,S_i,S_i')$ et $n:z\rightarrow y$ une image
inverse de $y$. Les restrictions $(mn)^*$ et $n^*m^*$ coincident
sur les objets et on suppose qu'elles coincident modulo un
morphisme de $F_{i+2}(S_1,S_1',..,U,z)\rightarrow
E_{i+2}(S_1,S_1',..,U,z)$.

On notera souvent de mani\`ere abusive $m^*$ et $m^*(U,V)$ par
$m^*$.

\medskip

On notera la $n+1-$tour gerb\'ee par $E(L_1,...,L_n,L_{n+1})$ et
on dira que $E_1$ est sa base.

\bigskip

\bigskip

On remarque que pour construire $E(L_1,...,L_n)$, on a construit
r\'ecursivement une famille de tours gerb\'ees
$(E(L_1)$,...,$E(L_1,..,L_i)$,...,$E(L_1,..,L_n))$. On appellera
cette famille de tours gerb\'ees la famille associ\'ee
$E(L_1,...,L_n)$.

\medskip

\medskip

{\bf D\'efinition 2.3}

Une tour gerb\'ee infinie de base $E_1$ est la donn\'ee pour tout
entier $n\geq 1$, d'une $n-$tour gerb\'ee $E(L_1,...,L_n)$ telle
que la famille associ\'ee \`a la $n+1-$tour gerb\'ee
$E(L_1,...,L_{n+1})$ est $(E(L_1),...,E(L_1,...,L_{n+1}))$.

\bigskip

{\bf Remarque.}

On supposera que les cat\'egories $F_{i+1}(S_1,S_1',...,S_i,S_i')$
et $E_{i+1}(S_1,S_1'..,S_i,S_i')$ sont $U-$petites, pour tout $i$
o\`u $U$ est un univers donn\'e. Ceci implique l'existence de la
r\'eunion des $F_{i+1}(S_1,S_1',..,S_i,S_i')$ et de celle des
$E_{i+1}(S_1,S_1'..,S_i,S_i')$.

\medskip

{\bf D\'efinition 2.4.}

Soient $E(L_1,...,L_n)$ une $n-$tour gerb\'ee, $S_1$ un objet de
$E_1$,  et $S_1'$ un objet de ${F_1}_{S_1}$ on peut d\'efinir une
famille de $n-1-$tours gerb\'ees   de la mani\`ere suivante:

On  pose $ E_{1S_1S_1'}(L_2)= F_2(S_1,S_1')\rightarrow
E_2(S_1,S_1')$,  $E_{iS_1S_1'}(L_2,..,L_{i+1})$, est la r\'eunion
des gerbes $F_{i+1}(S_1,S_1',...,S_i,S_i')\rightarrow
E_{i+1}(S_1,S_1',..,S_i,S_i')$.

La famille des $n-1-$tours gerb\'ees $E_{S_1S_1'}(L_2,...,L_n)$
est appel\'ee la fibre au-dessus de $S_1$.

\medskip

{\bf D\'efinition 2.5.}

Soit $E(L_1,...,L_n)$ une $n-$tour gerb\'ee  de base $E_1$, on
dira que $E'(L_1,...,L_n)$ est une sous $n-$tour gerb\'ee de
$E(L_1,...,L_n)$  si et seulement si $E'(L_1,...,L_n)$ est une
$n-$tour gerb\'ee au-dessus de $E_1$ telle que $F_1'\rightarrow
E_1$ est une sous-gerbe de $F_1\rightarrow E_1$, pour tout objet
$S_1$ de $E_1$ et $S_1'$ de ${F_1}'_{S_1}$,
$F_2'(S_1,S_1')\rightarrow E_2'(S_1,S_1)$ est une sous-gerbe de
$F_2(S_1,S_1')\rightarrow E_2(S_1,S_1'),...,
F_i'(S_1,S_1',..,S_{i-1},S_{i-1}')\rightarrow
E_i'(S_1,S_1'..,S_{i-1},S_{i-1}')$ est une sous-gerbe de
$F_i(S_1,S_1',..,S_{i-1},S_{i-1}')\rightarrow
E_i(S_1,S_1'..,S_{i-1},S_{i-1}')$,...,
$F_n'(S_1,S_1',..,S_{n-1},S_{n-1}')\rightarrow
E_n'(S_1,S_1'..,S_{n-1},S_{n-1}')$ est une sous-gerbe de
$F_n(S_1,S_1',..,S_{n-1},S_{n-1}')\rightarrow
E_n(S_1,S_1'..,S_{n-1},S_{n-1}').$

\medskip

{\bf D\'efinition 2.6.}

Un morphisme entre deux tours gerb\'ees  $E(L_1,...,L_n)$ et
$E'(L_1',...,L_n')$  de bases $E_1$ et $E_1'$ est la donn\'ee

d'un morphisme $f_1$ entre les gerbes $F_1\rightarrow E_1$ et
$F_1'\rightarrow E_1'$,

pour tout objets $S_1$ de $E_1$, et $S_1'$ de $F_{S_1}$ d'un
morphisme $f(S_1,S_1')$ entre les gerbes $F_2(S_1,S_1')\rightarrow
E_2(S_1,S_1')$ et $F_2'(f(S_1),f(S_1'))\rightarrow
E_2'(f(S_1),f(S_1'))$,...,

 pour tout objet $S_{i+1}$ de
$E_{i+1}(S_1,S_1'...,S_i,S_i')$ et $S_{i+1}'$ de
$F_{i+1}(S_1,S_1',..,S_i,S_i')_{S_{i+1}}$ d'un morphisme
$f(S_1,S_1'...,S_{i+1},S_{i+1}')$ entre les gerbes
$F_{i+2}(S_1,S_1',..,S_{i+1},S_{i+1}')\rightarrow
E_{i+2}(S_1,S_1'..S_{i+1},S_{i+1}')$ et
$F'_{i+2}(f(S_1),..,f(S_1,S_1'..,S_i,S_i')(S_{i+1}),f(S_1,S_1',..,S_i,S_i')
(S_{i+1}'))\rightarrow
E'_{i+2}(f(S_1),f(S_1'),..,f(S_1,S_1'..,S_i,S_i')(S_{i+1}),f(S_1,S_1,..,S_i,S_i')
(S_{i+1}'))$.

On commettra souvent l'abus de notation de noter
$f(S_1,S_1',..,S_i,S_i')$ par $f(S_i')$.

Soient $S_{i+1}$ un objet de $E_{i+1}(S_1,S_1',..,S_i,S_i')$,
$T_1$, $T_2$ des objets de
$F_{i+1}(S_1,S_1',..,S_i,S_i')_{S_{i+1}}$ et $u:T_1\rightarrow
T_2$ un morphisme;
 On suppose aussi la propri\'et\'e suivante v\'erifi\'ee:
$$
f(S_i')(u^*)=(f(S_i')(u))^*
$$
On suppose aussi que les morphismes commutent avec les
restrictions.

\bigskip

Soient $f:E(L_1,...,L_n)\rightarrow E'(L'_1,...,L'_n)$ un
morphisme de tours gerb\'ees, $S_i$ un objet de
$E_i(S_1,S_1',..,S_{i-1},S_{i-1}')$, $S_i'$ un objet de
$F_i(S_1,S_1',..,S_{i-1},S_{i-1}')_{S_i}$ et $l_i$ un \'el\'ement
de $L_i$.

Le morphisme $f$ d\'efinit un morphisme de liens
$f(S_i,L_i):L_i\rightarrow L_i$ entre $L_i$ et $L_i'$  par
$$
f(S_{i-1}')(l_iS_i')=f(S_i,L_i)(l_i)f(S_{i-1}')(S_i').
$$

La notion de morphisme entre tours gerb\'ees permet de d\'efinir
une notion d'\'equivalence  entre tours gerb\'ees.

\medskip

{\bf D\'efinition 2.7.}

On dira que deux $n-$tours gerb\'ees
 $E(L_1,...,L_n)$ et $E'(L_1,...,L_n)$ de base $E_1$ sont \'equivalentes
s'il existe un morphisme de tour gerb\'ee entre $E(L_1,...,L_n)$
et $E'(L_1,...,L_n)$ tel que le morphisme  $f(S_1,...,S_i')$ est
une \'equivalence de gerbes.

On note $H^{n+1}(E_1,L_1,...,L_n)$ l'ensemble des classes
d'\'equivalence des $n-$tours gerb\'ees de base $E_1$ et de liens
associ\'ees $L_1$,...,$L_n$.

\medskip

{\bf D\'efinition 2.8.}

La classe d'\'equivalence d'une  tour gerb\'ee $E(L_1,...,L_n)$
au-dessus de $E_1$ est triviale  (ou $n$ triviale) si et seulement
si il existe une sous tour gerb\'ee $E'(L_1,...,L_n)$ de
$E(L_1,...,L_n)$ telle que  les fibres de
$F_n'(S_1,S_1',..,S_{n-1},S_{n-1}')\rightarrow
E_n'(S_1,S_1'..,S_{n-1},S_{n-1}')$ n'ont qu'un \'el\'ement.

On dira qu'une $n-$tour gerb\'ee $E(L_1,...,L_n)$ est
$n-i$triviale si et seulement si il existe une sous $n-$tour
gerb\'ee $E'(L_1,...,L_n)$ de $E(L_1,...,L_n)$ telle que pour
$j\geq i$, les fibres de  $F'_j(S_1',...,S_{j-1}')\rightarrow
E'_j(S_1',...,S_{j-1}')$  n'ont qu'un \'el\'ement.

\bigskip

{\bf Le cas commutatif.}

\medskip

Dans cette partie, on suppose que les liens $L_1$,...,$L_n$ sont
des faisceaux en groupes commutatifs. On va exprimer les espaces
de cohomologies $H^{n+1}(E,L_1,...,L_n)$ en termes de cocycles de
Cech usuels.

\medskip

On suppose que la topologie de $E_1$ est d\'efinie par une famille
couvrante $(X_i\rightarrow X)_{i\in I}$. On suppose aussi que la
cohomologie de Cech d'un faisceau $L$ sur $E_1$ d\'efinit par
$(X_i\rightarrow X)_{i\in I}$ est la cohomologie de Cech de $L$
sur $E_1$.

\medskip

Posons $X_{i_1...i_k}=X_{{i_1}}\times_X...\times_X X_{{i_k}}$.

Si $x_{i_1}$ est un objet de ${F_1}_{X_{i_1}}$, on notera
$x_{i_1}^{i_2..i_k}$ sa restriction \`a
${F_1}_{X_{i_1}\times_X..\times_X X_{i_k}}$.

Si $u_{i_1i_2}$ est une fl\`eche entre $x_{i_2}^{i_1}$ et
$x_{i_1}^{i_2}$, on notera $u_{i_1i_2}^{i_3..i_k}$ sa restriction
entre $x_{i_2}^{i_1..i_k}$ et $x_{i_1}^{i_2..i_k}$.

 Soit
$x_i$ un objet de ${F_1}_{X_i}$, on note $u_{i_1i_2}$ une fl\`eche
entre les restrictions respectives $x_{i_2}^{i_1}$ et
$x_{i_1}^{i_2}$ de $x_{i_2}$ et $x_{i_1}$ \`a
${F_1}_{X_{i_1i_2}}$, la fl\`eche
$c_{i_1i_2i_3}=u_{i_2i_3}^{i_1}{(u_{i_1i_3}^{i_2})}^{-1}u_{i_1i_2}^{i_3}$
repr\'esente un $2-$cocycle de Cech du faisceau $L_1$. Son bord
$d(c_{i_1i_2i_3})$ est donc nul.

 On note $x_{i_1i_2}'$ un objet de $E_2(X_{i_1i_2},x_{i_2}^{i_1})$ et $x_{i_1i_2}$ un
 objet de sa fibre,
 $(c_{i_1i_2i_3})^*$ est un morphisme de la gerbe
$F_2(X_{i_1i_2i_3},x_{i_2}^{i_1i_3})\rightarrow
E_2(X_{i_1i_2i_3},x_{i_2}^{i_1i_3})$ qui agit sur
$x_{i_1i_2}'^{i_3}$.

On a vu que les chaines $c_{i_2i_3i_4}^{i_1}$,
$c_{i_1i_3i_4}^{i_2}$, $c_{i_1i_2i_3}^{i_4}$ peuvent \^etre
consid\'erer comme \'el\'ements de $L_1$. On peut dont aussi les
consid\'erer leur image par $*$ comme morphisme de la gerbe
$F_2(X_{i_1i_2i_3i_4}, x_{i_2}^{i_1i_3i_4})\rightarrow
E_2(X_{i_1i_2i_3i_4}, x_{i_2}^{i_1i_3i_4})$.
 Via ces identifications, le bord de $(c_{i_1i_2i_3})^*$
  est un morphisme de la gerbe
 $F_2(X_{i_1i_2i_3i_4},..,x_{i_2}^{i_1i_3i_4})
 \rightarrow E_2(X_{i_1i_2i_3i_4},..,x_{i_2}^{i_1i_3i_4})$ qui
 pr\'eserve ${x_{i_1i_2}'}^{i_3i_4}$.
  Il
 s'identifie donc \`a une chaine $c_{i_1i_2i_3i_4}$ de $L_2$.

\medskip

{\bf Lemme 2.9.} {\it La chaine $c_{i_1i_2i_3i_4}$ est un
$3-$cocycle.}

\medskip

{\bf Preuve.}

Evaluons le bord de  $c_{i_1..i_4}$, on obtient:
$$
d(c_{i_1..i_4})=\sum_{j=1}^{j=5}(-1)^{j}c_{i_1..\hat i_j..i_5}
$$
$$
=\sum_{j=1}^{j=5}(-1)^j(\sum_{k=1}^{k=j-1}(-1)^kc_{i_1..\hat
i_k..\hat i_j..}^*+\sum_{k=j+1}^{k=5}(-1)^{k+1}c_{i_1..\hat
i_j\hat i_k..i_5}^*)=0.
$$

\medskip

Supposons d\'efini le $L_j$ $j+1-$cocycle $c_{i_1...i_{j+2}}$,
alors $(c_{i_1..i_{j+2}})^*$ est un morphisme de la gerbe
$F_{j+1}(X_{i_1..i_{j+2}},x_{i_2}^{1..j+2},..,
x_{i_2..i_{j}}'^{i_{j+1}i_{j+2}},
x_{i_2..i_{j}}^{i_{j+1}i_{j+2}})\rightarrow
E_{j+1}(X_{i_1..i_{j+2}},x_{i_2}^{2..j+2},..,
x_{i_2..i_{j}}'^{i_{j+1}i_{j+2}},
x_{i_2..i_{j}}^{i_{j+1}i_{j+2}})$. O\`u $x_{i_1..i_{j}}'$ un objet
de
$E_j(X_{i_1..i_{j+2}},..,x_{i_2..i_{j-1}}^{i_ji_{j+1}i_{j+2}})$,
et $x_{i_1..i_{j}}$ un objet de sa fibre.

Si $i_{j+1}=i_{j+2}$, alors peut supposer que $c_{i_1...i_{j+2}}$
est nul.

Chaque fl\`eche
$c_{i_2..i_{j+3}}^{i_1}$,...,$c_{i_1..i_{j+2}}^{i_{j+3}}$,
s'identifie \`a un \'el\'ement de $L_j$.

 Le bord
 $d(c_{i_1..i_{j+2}})^*=c_{i_1..i_{j+3}}$ s'identifie donc \`a
  un morphisme de la
 gerbe
 $F_{j+1}(X_{i_1..i_{j+3}},x_{i_2}^{2..j+2},..,
x_{i_2..i_{j}}^{i_{j+1}i_{j+2}i_{j+3}}) \rightarrow
E_{j+1}(X_{i_1..i_{j+3}},x_{i_2}^{2..j+2},..,
x_{i_2..i_{j}}^{i_{j+1}i_{j+2}i_{j+3}})$ qui pr\'eserve un objet
$x_{i_1..i_{j+1}}'$ de sa base.

Il s'identifie donc \`a un \'el\'ement de  $L_{j+1}$.

\medskip

{\bf Lemme 2.10.}

{\it La famille $c_{i_1...i_{j+3}}$ est un $j+2$ cocycle.}

\medskip

{\bf Preuve.}

Evaluons le bord de $c_{i_1...i_{j+3}}$, on a
$$
d(c_{i_1...i_{j+3}})=\sum_{k=1}^{k=j+4}(-1)^k c_{i_1..\hat
i_k..i_{j+4}}
$$
=
$$
\sum_{k=1}^{k=j+4}(-1)^k(\sum_{l=1}^{l=k-1}(-1)^lc_{i_1..\hat
i_l..\hat i_k..}^*+\sum_{l=k+1}^{l=j+4}(-1)^{l+1}c_{i_1..\hat
i_k\hat i_l..i_{j+4}}^*)=0.
$$

On d\'efinit ainsi r\'ecursivement un $n+1-$cocycle
$c_{i_1...i_{n+2}}$.

\medskip

{\bf Remarque.}

Le $j+1-$cocycle $c_{i_1..i_{j+2}}$ est d\'efini modulo un
\'el\'ement de $L_j$. Mais sa classe de cohomologie ne d\'epend
pas du repr\'esentant choisi.

Soit $(E_n(L_1,...,L_n))_{n\in {\N}}$ une tour gerb\'ee infinie
telle que chaque faisceau $L_n$ est ab\'elien. Notons $c_{n+1}$ le
cocycle associ\'e \`a la $n-$tour gerb\'ee $E(L_1,...,L_n)$. On
note $c$  la famille $(c_n)_{n\in {\N}}$.

\bigskip

{\bf Proposition 2.11.}

{\it Soit $E'(L_1,...,L_n)$ une sous $n-$tour gerb\'ee de
$E(L_1,...,L_n)$, alors les classes de cohomologies des familles
de cocycles $(c'_2,...,c'_{n+1})$ et $(c_2,...,c_{n+1})$ qui
d\'ecrivent respectivement les tours gerb\'ees  $E'(L_1,...,L_n)$
et $E(L_1,...,L_n)$ coincident.}

\medskip

{\bf Preuve.}

Soit $(X_i\rightarrow X)_{i\in I}$ une famille couvrante de la
topologie de $E_1$, $c_2$ est construite en choisissant dans
$F_{1_{X_i}}$ un objet $x_i$ qu'on peut supposer appartenir \`a
$F'_{1_{X_i}}$, il en r\'esulte que $c_2$ et $c_2'$ coincident.

Supposons que $c_{j+1}$ coincide avec $c'_{j+1}$, pour construire
$c_{j+2}$, on identifie $(c_{i_1...i_{j+2}})^*$ \`a un morphisme
de la gerbe $F_{j+1}(X_{i_1..i_{j+3}},x_{i_2}^{2..j+3},..,
x_{i_2..i_{j}}^{i_{j+1}i_{j+2}i_{j+3}}) \rightarrow
E_{j+1}(X_{i_1..i_{j+3}},x_{i_2}^{2..j+3},..,
x_{i_2..i_{j}}^{i_{j+1}i_{j+2}i_{j+3}})$. On peut aussi supposer
r\'ecursivement que c'est un morphisme de la gerbe
$F'_{j+1}(X_{i_1..i_{j+3}},x_{i_2}^{2..j+3},..,
x_{i_2..i_{j}}^{i_{j+1}i_{j+2}i_{j+3}}) \rightarrow
E'_{j+1}(X_{i_1..i_{j+3}},x_{i_2}^{2..j+3},..,
x_{i_2..i_{j}}^{i_{j+1}i_{j+2}i_{j+3}})$. Ceci implique que
$c_{j+2}=c'_{j+2}$.

\medskip

{\bf Corollaire 2.12.}

{\it Supposons que la $n-$gerbe $E(L_1,..,L_n)$ soit triviale,
alors la classe de cohomologie du $n+1-$cocycle de Cech qu'on
vient de d\'efinir est triviale.}

\medskip

{\bf Preuve.}

Supposons que la $n-$gerbe $E(L_1,...,L_n)$ soit triviale, Il
existe alors une sous $n-$gerbe $E'(L_1,...,L_n)$ de
$E(L_1,...,L_n)$ telle que la fibre $F'_n(S_1,...,S_n)\rightarrow
E'_n(L_1,...,L_n)$ a un \'el\'ement. Il suffit de montrer que le
cocycle $c'_{n+1}$ est trivial. Dans ce cas $c'_{n+1}$ est
construit en identifiant le bord de $c'_{n}$ a une famille de
morphismes de gerbes triviales. Le r\'esultat se d\'eduit en
remarquant que les $c_{i_1..i_{n+2}}$, pour $i_1=...=i_{n-1}$ sont
l'obstruction \`a la trivialit\'e des gerbes
$F_n(S_1,...,S_{n-1}')\rightarrow E_n(S_1,...,S_{n-1}')$.

\medskip

{\bf Th\'eor\`eme 2.13.}

{\it Soit $E_1$ une cat\'egorie et $L_1$,...,$L_n$ des faisceaux
de $E_1$ ab\'eliens de $E_1$.
 L'ensemble des classes d'isomorphismes des $n-$tours
gerb\'ees ab\'eliennes s'identifie au groupe de Cech
$H^{n+1}(E_1,L_n)$ par l'application qui \`a une $n-$tour gerb\'ee
$E(L_1,...,L_n)$ associe la classe de cohomologie  $[c_{n+1}]$ de
$c_{n+1}$.}

\medskip

{\bf Preuve.}

Montrons d'abord que si  $[c_{n+1}]$ est nulle alors la tour
gerb\'ee est triviale, ou de mani\`ere \'equivalente que les
gerbes $F_n(S_1,..,S_{n-1})\rightarrow E_n(S_1,..,S_{n-1})$ sont
triviales.

Le cocycle $c_{n+1}$ est consid\'er\'e comme famille de morphismes
$c_{i_1..i_{n+2}}$ des gerbes
$F_n(X_{i_1..i_{n+2}},..,x_{i_2..i_{n-1}}^{i_{n+1}i_{n+2}})\rightarrow
E_n(X_{i_1..i_{n+2}},..,x_{i_2..i_{n-1}}^{i_{n+1}i_{n+2}})$.

On choisit $i_1=..=i_{n-1}=i_0$. Dans ce cas, le cocycle
$c_{i_1..i_{n+2}}$ se r\'eduit \`a l'obstruction \`a la
trivialit\'e de la gerbe
$F_n(X_{i_1..i_{n+2}},..,x_{i_2..i_{n-1}}^{i_{n}i_{n+1}i_{n+2}})\rightarrow
E_n(X_{i_1..i_{n+2}},..,x_{i_2..i_{n-1}}^{i_{n}i_{n+1}i_{n+2}})$.

Il faut maintenant montrer qu'\`a partir de toute classe
$c_{n+1}$, on peut construire une tour gerb\'ee.

Pour cela, il faut d'abord montrer que l'ensemble des classes
d'\'equivalences des $L_1$ gerbes est en bijection avec
$H^2(E_1,L_1)$, puis par r\'ecurence que l'ensemble des classes
d'\'equivalences des $i-$tours gerb\'ees $E(L_1,..,L_i)$ de base
$E_1$ est en bijection avec $H^{i+1}(E_1,L_i)$.

Montrons donc que l'ensemble des classes d\'equivalences des
gerbes de base $E_1$ est en bijection avec $H^2(E_1,L_1)$.

Soit $c_{i_1i_2i_3}$ un $2-$cocycle repr\'esentant une classe de
$H^2(E_1,L_1)$.

Consid\'erons une famille de gerbes triviales $P_i\rightarrow X_i$
de lien $L_1$  telle qu'on ait des morphismes $c_{ij}$ entre les
restrictions respectives de $P_j$ \`a $X_{ij}$ et de $P_i$ \`a
$X_{ij}$ de sorte que
$$
c_{i_2i_3}c_{i_1i_3}^{-1}c_{i_1i_2}=c_{i_1i_2i_3}.
$$
Pour s'assurer de l'existence de $c_{i_1i_2}$, il suffit de
remarquer qu'\'etant donn\'e $X_{i_0}$, on peut d\'efinir la gerbe
$P_{i_0}$ sur la cat\'egorie au-dessus de $X_{i_0}$ de sorte que
ses restrictions $P_{i_0i_1}$ au-dessus de $X_{i_1}\times_X
X_{i_0} \rightarrow {X_{i_0}}$ soit la gerbe triviale de lien
$L_1$. On pose alors $u_{i_1i_2}=c_{i_1i_2i_0}$ comme morphisme
entre les restrictions respectives de $P_{i_2i_0}$ et $P_{i_1i_0}$
\`a $X_{i_1}\times_X X_{i_2}\times_XX_{i_0}\rightarrow X_{i_0}$.

\medskip

Supposons que pour tout \'el\'ement $[c_{i+1}]$ de
$H^{i+1}(E_1,L_{i})$, il existe une tour gerb\'ee $E(L_1,...,L_i)$
de base $E_1$ et de cocycle classifiant $c_{i+1}$.

Soit $c_{i+2}$ un \'el\'ement de $H^{i+2}(E_1,L_{i+1})$. On va
construire \`a partir d'une $i-$tour gerb\'ee quelconque
$E(L_1,...,L_i)$, une tour $i+1-$tour gerb\'ee de cocycle
classifiant $c_{i+2}$.

Le cocycle classifiant de $E(L_1,...,L_i)$ est donn\'ee par une
famille de chaines $c_{i_1...i_{j+2}}$ qui peuvent \^etre
consid\'er\'ees comme morphismes de gerbes
$F_j(X_{i_1..i_{j+2}},x_{i_2}^{2..j+2},..,x_{i_2..i_{j-1}}^{i_{j+1}i_{j+2}})\rightarrow
E_j(X_{i_1..i_{j+2}},x_{i_2}^{2..j+2},..,x_{i_2..i_{j-1}}^{i_{j+1}i_{j+2}})$.

A tout objet $x_{i_1...i_{j}}'$ de
$E_j(X_{i_1..i_{j+2}},x_{i_2}^{2..j+2},..,x_{i_2..i_{j-1}}^{i_{j+1}i_{j+2}})$
et $x_{i_1..i_{j}}$ de sa fibre, on associe une gerbe triviale
$F_{j+1}(X_{i_1..i_{j+2}},x_{i_2}^{2..j+2},..,
x_{i_1..i_{j}}',x_{i_1..i_{j}}) \rightarrow
E_{j+1}(X_{i_1..i_{j+2}},x_{i_2}^{2..j+2},..,x_{i_2..i_{j-1}}^{i_{j+1}i_{j+2}},
x_{i_1..i_{j}}',x_{i_1..i_{j}})$.

On d\'efinit des morphismes
$h_{i_1...i_{j+2}}=(c_{i_1..i_{j+2}})^*$ de cette gerbe
 de sorte que
$d(h_{i_1..i_{j+2}})=c_{i_1..i_{j+3}}$.

Pour s'assurer de l'existence des applications $h_{i_1..i_{j+2}}$,
on les d\'efinit pour $i_0$ fix\'e sur  leurs restrictions
$F_{j+1}(X_{i_1i_2..i_{j+2}i_0},x_{i_2}^{2..j+3},..,x_{i_2..i_{j}}'^{i_{j+3}},
,x_{i_1..i_{j}}^{i_{j+3}}) \rightarrow
E_{j+1}(X_{i_1i_2..i_{j+2}i_0},x_{i_2}^{2..j+3},..,x_{i_2..i_{j}}^{i_{j+3}})$
 par
$h_{i_2..i_{j+3}}^{i_0}=c_{i_1i_2..i_{j+2}i_0}$.

\bigskip

{\bf 3. Suites spectrales et tours gerb\'ees.}

\medskip

Le but de cette partie est d'appliquer la m\'ethode des suites
spectrales pour \'etudier les groupes de cohomologie des  tours
gerb\'ees ab\'eliennes.

Soit $E(L_0,...,L_n,..)_{n\in{\N}}$ une $\infty-$tour gerb\'ee
o\`u la famille $(L_n)_{n\in {\N}}$ est une famille de faisceaux
ab\'eliens au-dessus de la cat\'egorie $E_1$. On suppose que la
topologie de $E_1$ est d\'efinie par la famille couvrante
$(X\rightarrow X_i)_{i\in I}$.

Posons $L=\oplus_{i\geq 0}L_i$, et notons $(C^*(X,L),d)$ le
complexe des chaines de Cech d\'efinies sur $E_1$ et \`a valeurs
dans $L$.

On peut munir ce complexe   des chaines de Cech  de la filtration
suivante:

$$
K_p=C^*(X,\oplus_{q\geq p} L_q)
$$
et de la graduation
$$
K^p= C^p(X,L)
$$
On va calculer les termes de la suite spectrale associ\'ee \`a
cette filtration et \`a cette graduation.

 Notons $ Z^p_r=\{x\in K_p: d(x)\in K_{p+r}\},$ $B^p_r=d(K_{p-r}\cap K_p),$
   et $E^p_r= {Z^p_r\over
{Z^{p+1}_{r-1}+B^p_{r-1}}}$.

On suppose dans la suite que $r\geq 1$.

\medskip

D\'eterminons $Z^p_r$.

Soit $x$ un \'el\'ement de $K_p$, $x$ appartient \`a $Z^p_r$ si et
seulement si $d(x)$ appartient \`a $K_{p+r}$. Posons
$x=x_{i_0}+...+x_{i_n}$, o\`u $x_{i_l}$ est la composante
homog\`ene de $x$ \`a valeurs dans $L_l$. On a $d(x_l)$ est un
\'el\'ement de $K_{p+r}$ si et seulement si $l\geq p+r$ ou bien si
$d(x_l)=0$. On en d\'eduit que $x$ appartient \`a $Z^p_r$ si et
seulement si ses composantes $x_{i_j}$ telles que $j<p+r$ sont des
cocycles.

\medskip

D\'eterminons $B^p_r$.

Soit $x$ un \'el\'ement de $K_{p-r}$, une de ses composantes
homog\`enes $x_l$ \`a valeurs dans $L_l$ appartient \`a $B^p_r$ si
et seulement si $d(x_l)$ appartient \`a $K_p$ ceci est
\'equivalent \`a dire que $d(x_l)$ est nul ou bien que $l\geq
p+r$. Il en r\'esulte que $B^p_r= d(K_{p})$.

\medskip

D\'eterminons $E^p_r$.

On a $Z^p_r=Z^{p+1}_{r-1}\oplus Z(K)\cap C(X, L_p).$ On en
d\'eduit que $E^p_r= H(X,L_p)$.

\bigskip

Posons maintenant $Z^{pq}_r=Z^p_r\cap K^{p+q}$,
$B_r^{pq}=B_r^p\cap K^{p+q}$ et $E^{pq}_r={Z^{pq}_r\over
B^{pq}_{r-1}+Z^{p+1,q-1}_{r-1}}$.

D\'eterminons $Z^{pq}_r$.

Soit $x$ un \'el\'ement de $Z^p_r$, une de ses composantes
homog\`enes $x_l$ \`a valeurs dans $L_l$   appartient \`a
$Z_r^{pq}$ si et seulement si $x_l$ est une $p+q-$chaine \`a
valeurs dans $K_p$,  et $d(x_l)$ appartient \`a $K_{p+r}$. On a vu
que $d(x_l)$ appartient \`a $K_{p+r}$ si et seulement si $x_l$ est
un cocycle ou bien $l\geq p+r$. On en d\'eduit que
$Z^{pq}_r=C^{p+q}(X,L_{p+r})\oplus Z^{p+q}(X,L_p\oplus..\oplus
L_{p-1})$.

\medskip

D\'eterminons maintenant $B^{pq}_r$.

Soit $x$ un \'el\'ement de $B^{pq}_r$, une de ses composantes
$x_l$ \`a valeurs  dans $L_l$ appartient \`a $B^{pq}_r$, si et
seulement si c'est une $p+q-$chaine   et il existe $y$ appartenant
\`a $K_{p-r}$ tel que $d(y_l)=x_l$. On en d\'eduit que
$B^{pq}_r=d(C^{p+q-1}(X,K_p))$.

\medskip

D\'eterminons $E^{pq}_r$.

L'espace vectoriel $Z^{pq}_r$ est la somme de $Z^{p+1,q-1}_{r-1}$
et de $Z^{p+q}(X,L_p)$. On en d\'eduit que
$E^{pq}_r=H^{p+q}(X,L_p)$.

\medskip

On va noter $Z^p_{\infty}$ l'ensemble des cycles contenus dans
$K_p$, $B^p_{\infty}$ l'ensemble des bords contenus dans $K_p$ et
enfin
$E^p_{\infty}={Z^p_{\infty}\over{Z^{p+1}_{\infty}+B^p_{\infty}}}$.

On remarque que $E^p_{\infty}=H(X,L_p)$.

\medskip

La proposition suivante se d\'eduit de [God] p. 84 Th\'eor\`eme
4.6.1.

\bigskip

{\bf Proposition 3.1.}

 {\it Supposons maintenant qu'il existe un
entier $n\geq r$ tel que $H^{p+q}(X,L_p)=E^{pq}_r=0$ pour $p\neq
0,n$ et un entier $s$ tel que $L_n=0$ si $n>s$.

On obtient alors la suite exacte suivante:
$$
...\rightarrow H^{n}(X,L_n)\rightarrow H^i(X,L)\rightarrow
H^i(X,L_0)\rightarrow H^{i+1}(X,L_n)\rightarrow
H^{i+1}(X,L)\rightarrow...
$$
}

\bigskip

{\bf 4. Application des tours gerb\'ees \`a la g\'eom\'etrie
affine.}

\bigskip

Dans cette partie, on va appliquer la notion de tour gerb\'ee au
probl\`eme qui a motiv\'e sa construction: la classification des
s\'eries de compositions des vari\'et\'es affines.

\medskip

Soient $(M_1,\nabla_{M_1})$,
$(F_1,\nabla_{F_1})$,...,$(F_{n-1},\nabla_{F_{n-1}})$ des
vari\'et\'es affines compactes et compl\`etes, on se propose de
classifier les s\'eries de compositions de la forme
$$
(M_n,\nabla_{M_n})\rightarrow
(M_{n-1},\nabla_{M_{n-1}})...\rightarrow
(M_2,\nabla_{M_2})\rightarrow (M_1,\nabla_{M_1})
$$
o\`u pour tout $i$ l'application
$$
f_i:(M_{i+1},\nabla_{M_{i+1}})\rightarrow (M_i,\nabla_{M_i})
$$
est une fibration affine c'est \`a dire une application affine
surjective, dont l'espace source une vari\'et\'e affine  (compacte
et compl\`ete)  et dont la fibre  est diff\'eomorphe \`a $F_i$, et
est munie d'une structure affine d'holonomie lin\'eaire celle de
$(F_i,\nabla_{F_i})$.

On notera $m_1$ la dimension de $M_1$, $l_i$ la dimension de
$F_i$, et $H^1(\pi_1(F_i),{\R}^{l_i})$ le premier groupe de
cohomologie  de $\pi_1(F_i)$ relatif \`a l'holonomie lin\'eaire de
$(F_i,\nabla_{F_i})$.

L'application $f_{i+1}$ induit une repr\'entation
$:\pi_1(F_i)\rightarrow Gl(H^1(\pi_1(F_{i+1}),{\R}^{l_{i+1}}))$
(Voir [T6]).

Les donn\'ees \`a partir desquelles se feront la classification
seront:

- La vari\'et\'e affine $(M_1,\nabla_{M_1})=(F_0,\nabla_{F_0})$,

- Une repr\'esentation $\pi_i:\pi_1(F_i)\rightarrow
Aff({\R}^{m_1+l_1..+l_i})$ dont la partie lin\'eaire pr\'eserve
$0\times {\R}^{l_i}$ et qui se projette sur ${\R}^{m_1+..l_{i-1}}$
en l'identit\'e. On suppose que le quotient de
${\R}^{m_1+l_1+..+l_i}$ par $\pi_1(F_i)$ est une fibration
au-dessus de ${\R}^{m_1+l_1+..+l_{i-1}}$ dont les fibres sont
diff\'eomorphes \`a $F_i$.

La partie lin\'eaire de la restriction de $\pi_i$ \`a $0\times
{\R}^{l_i}$ est l'holonomie lin\'eaire des structures affines de
$F_i$ consid\'er\'ees.

- Soit $N(I(\pi_{i+1}))$ le sous-groupe du normalisateur de
l'image de $\pi_{i+1}$ dans $Aff({\R}^{m_1+l_1+..+l_{i+1}})$ dont
la partie lin\'eaire d'un de ses \'el\'ements quelconque $g$
pr\'eserve $0\times {\R}^{l_{i+1}}$ et tel que l'action de $g$
induite sur ${\R}^{m_1+..l_i}$ est l'identit\'e. On se donne aussi
des repr\'esentations
$$
\pi_i':\pi_1(F_i)\rightarrow N(I(\pi_{i+1}))
$$
celle ci induit une repr\'esentation
$$
\pi_i'':\pi_1(F_i)\rightarrow
Gl(H^1(\pi_1(F_{i+1}),{\R}^{l_{i+1}}))
$$

o\`u   $Gl(H^1(\pi_1(F_i),{\R}^{l_i}))$ le groupe des
automorphismes lin\'eaires de $H^1(\pi_1(F_i),{\R}^{l_i})$.

 La
repr\'esentation $\pi_i$ induit aussi sur $F_i$ un fibr\'e plat de
fibre type $H^1(\pi_1(F_{i+1}),{\R}^{l_{i+1}})$ qu'on note aussi
$\pi_i$.

\medskip

Soit $A_i=(Aff({\R}^{m_1+..+l_i})/\pi_1(F_i))_0$ la composante
connexe du groupe des automorphismes affines de
${\R}^{m_1+l_1+..+l_i}/\pi_1(F_i)$ qui pr\'eservent ses fibres et
se projettent sur ${\R}^{m_1+..l_{i-1}}$ en l'identit\'e, $\pi_i'$
d\'efinit au-dessus de $F_i$ un fibr\'e plat de fibre type
$A_{i+1}$, on note $h_{i+1i}$ le faisceau des sections affines de
ce fibr\'e. (section affines veut dire qu'elles se rel\`event
localement sur $N(I(\pi_{i+1}))$ en des sections de $F_i$ dont la
partie lin\'eaire ne d\'epend pas de la variable et la partie
translation en d\'epend de mani\`ere affine.) C'est un fibr\'e
au-dessus de $F_i$.

La repr\'esentation $\pi'_{i-1}$ d\'efinit un fibr\'e plat
au-dessus de $F_{i-1}$ de fibre type $h_{i+1i}$, qu'on note
$h_{i+1i-1}$ le faisceau des sections affines de ce fibr\'e.

R\'ecursivement, on peut d\'efinir les faisceaux $h_{i+11}$.

\bigskip

L'outil n\'ecessaire pour \'effectuer la classification est la
notion de $n-$tour gerb\'ee qu'on vient de d\'efinir. On va
associer \`a la famille $(M_1,\nabla_{M_1})$,
$(F_1,\nabla_{F_1})$,...,$(F_{n-1},\nabla_{F_{n-1}})$ la
$n-1-$tour gerb\'ee suivante:

Soit $U$ un ouvert de $M_1$, on note $C_1(U)$ la cat\'egorie dont
les objets sont les fibr\'es affines de base $U$, et tels que le
fibr\'e plat au-dessus de $U$, de base
$H^1(\pi_1(F_1),{\R}^{l_1})$ induit par un objet de $C_1(U)$ est
la restriction de $\pi_1''$ \`a $U$. Les morphismes des objets de
$C_1(U)$ sont les morphismes de fibr\'es affines c'est \`a dire
les isomorphismes affines entre espaces totaux qui respectent les
fibres et se projettent sur $U$ en l'identit\'e.

 Le faisceau de cat\'egories d\'efini sur $M_1$ par
$$
U\longrightarrow C_1(U)
$$
est une gerbe de lien $A_1$.

\medskip

Supposons $U$ simplement connexe.
 Soit $S_1$ un objet de $C_1(U)$, c'est un fibr\'e affine de base
 $U$. Sa fibre est diff\'eomorphe \`a $F_1$.
 Consid\'erons un ouvert $U(U_1)$  de $S_1$ diff\'eomorphe \`a
 $U\times U_1$ o\`u $U_1$ est un ouvert de $F_1$.
  On peut d\'efinir la cat\'egorie $C_2(U(U_1))$ dont les
 objets sont les fibr\'es affines au-dessus de $U(U_1)$ de
  fibre diff\'eomorphe \`a $F_2$,  tels que
 le fibr\'e plat de fibre type $H^1(\pi_1(F_2),{\R}^{l_2})$
 associ\'e soit induit par $\pi_2$.

  On vient de d\'efinir ainsi une $2-$tour gerb\'ee
  $C_2(A_1, A_2)$ au-dessus de
  $M_1$.

  Supposons d\'efinie la $i-$tour $E_i$ gerb\'ee, elle
  consiste en la donn\'ee d'une famille de fibr\'es affines de
  fibres $F_i$ au-dessus d'une famille d'ouvert
  $U(U_1,...,U_{i-1})$ diff\'eomorphe \`a
  $U\times U_1...\times U_{i-1}$, o\`u $U_k$ est un ouvert  de
  $F_k$, qui est simplement connexe si $k<i-1$.
  Soit $U_i$ un ouvert de $F_i$, ayant d\'efini $C_i$, on peut d\'efinir
  $C_{i+1}$ en associant \`a
  $U(U_1,..,U_i)$ ouvert d'un objet de $C_i(U(U_1,...,U_{i-1}))$
   diff\'eomorphe \`a $U( U_1,..,U_{i-1})\times U_i$
   (o\`u cette
  fois $U_{i-1}$ est simplement connexe) la
  cat\'egorie, $C_{i+1}(U( U_1,... U_i))$ dont les objets
  sont les fibr\'e affines au-dessus de $U(U_1,... U_i)$ de
  fibre type $F_{i+1}$ tel que le fibr\'e plat de fibre type
  $H^1(\pi_1(F_{i+1}),{\R}^{l_{i+1}})$ associ\'e soit induit par $\pi_i''$.

  Le faisceau de cat\'egories
  $$
  U(U_1,..., U_i)\longrightarrow C_{i+1}(U(U_1,..., U_i))
  $$
  est une gerbe au-dessus de $U(U_1,..., U_{i-1})\times F_i$
  de liens $L_{i+1}$.

\medskip

On d\'efinit  ainsi par r\'ecurrence une tour gerb\'ee
$E(L_1,...,L_{n-1})$.

\medskip

\medskip

{\bf Proposition 4.1.}

{\it Il existe une s\'erie de composition de
$(M_n,\nabla_{M_n})\rightarrow...\rightarrow (M_1,\nabla_{M_1})$
telle que la fibration $f_i:(M_{i+1},\nabla_{M_{i+1}})\rightarrow
(M_i,\nabla_{M_i})$ de fibre diff\'eomorphe \`a $F_i$ associ\'e au
donn\'ees d\'efinies en d\'ebut de ce num\'ero si
  et
seulement si la tour gerb\'ee associ\'ee est $n-1-$triviale.}

\medskip

{\bf D\'efinition. 4.2.}

Un fibr\'e affine $(M,\nabla_M)\rightarrow (B,\nabla_B)$ d'espace
total compact est dit affinement localement trivial, si et
seulement si l'holonomie d'une fibre ne d\'epend pas de celle-ci.

\medskip

En d'autres termes, ceci veut dire que le fibr\'e est un fibr\'e
diff\'erentiable localement trivial, et le cocycle de Cech qui
sert \`a le d\'efinir peut \^etre exprim\'e par des applications
affines de la fibre qu'on note $(F,\nabla_F)$.

Les s\'eries de compositions telles que $f_i$ est un fibr\'e
affine affinement localement trivial ont \'et\'e classifi\'ees
dans [T6]. Dans ce cas on peut se restreindre \`a des tours
gerb\'ees de liens commutatifs, et exprimer les invariants en jeu
par des cocycles de Cech usuels.

\bigskip

{\bf 5. Interpr\'etation de la cohomologie enti\`ere d'une
vari\'et\'e.}

\medskip

Les classes caract\'eristiques ont \'et\'e recemmemt utilis\'ees
par plusieurs math\'ematiciens afin d'\'etudier des objets
g\'eom\'etriques. Pour exploiter au mieux la th\'eorie des classes
caract\'eristiques, il faut donner une interpr\'etation
g\'eom\'etrique des groupes $H^n(M,{\Z})$, o\`u $M$ est une
vari\'et\'e diff\'erentiable. C'est le but de cette partie.

Soit $M$ une vari\'et\'e diff\'erentiable. On a vu que le groupe
$H^{n+1}(M,{\C}^*)$ s'identifie aux classes d'\'equivalences des
$n-$tours gerb\'ees $E_n({\C}^*_M,...,{\C}^*_M)$, o\`u
 ${\C}^*_M$, (resp ${\C}_M$) est le faisceau des fonctions
diff\'erentiables d\'efinies sur $M$ et \`a valeurs dans ${\C}^*$.
(resp. le faisceau des fonctions d\'efinies sur $M$ et \`a valeurs
dans ${\C}$). On a une suite exacte
$$
0\rightarrow {\Z}{\buildrel  i\over{\rightarrow}} {\C}_M{\buildrel
{exp } \over{\rightarrow}}{\C}^*{_M}\rightarrow 0.
$$
La suite exacte de cohomologie associ\'ee \`a cette suite exacte
de faisceaux donne un isomorphisme entre les groupes
$H^{n+1}(M,{\Z})$ et le groupe $H^n(M,{\C}^*_M)$. Il vient du fait
que $H^{n+2}(M,{\Z})$ est isomorphe au groupe des classes
d'\'equivalences des $n-$tours gerb\'ees
$E_n({\C}^*_M,...,{\C}^*_M)$ le r\'esultat suivant:

\medskip

{\bf Th\'eor\`eme 5.1.} {\it Soit $M$ une vari\'et\'e
diff\'erentiable, le groupe de cohomologie $H^{n+2}(M,{\Z})$ est
isomorphe aux classes d'\'equivalences des $n-$tours gerb\'ees
d\'efinies sur $M$.}

\bigskip

{\bf Acknowledgements.}

\medskip

The author would acknowledges the Abdus Salam ICTP, Trieste, Italy
for support.  This work has been done
 at ICTP.

\bigskip

\centerline{\bf Bibliography.}

\bigskip

[A] Atiyah, M. Complex analytic connections in fiber bundles.
Trans. Amer. Math. Soc 85 (1957)181-207.

[Bo] Borel, A. Linear algebraic groups. W.A. Benjamin, Inc, New
York-Amsterdam 1969.

[Bre] Breen, L. On the classification of $2-$gerbes and
$2-$stacks. Asterisque, 225 1994.

[Br] Bredon, G. E. Sheaf theory. McGraw-HillBook Co., 1967.

[Bry] Brylinski, J.L Loops spaces, Characteristic Classes and
Geometric Quantization, Progr. Math. 107, Birkhauser, 1993.

[Br-Mc] Brylinski, J.L, Mc Laughlin D.A, The geometry of degree
four characteristic classes and of line bundles on loop spaces I.
Duke Math. Journal. 75 (1994) 603-637.

[Ca] Carriere, Y. Autour de la conjecture de L. Markus sur les
vari\'et\'es affines. Invent. Math. 95 (1989) 615-628.

[De] Deligne, P. Theorie de Hodge III, Inst. Hautes Etudes Sci.
Publ. Math. 44 (1974), $5-77$.

[Du] Duskin, J. An outline of a theory of higher dimensional
descent, Bull. Soc. Math. Bel. S\'erie A 41 (1989) 249-277.

[F] Fried, D. Closed similarity affine manifolds. Comment. Math.
Helv. 55 (1980) 576-582.

[Gi] Giraud, J. Cohomologie non ab\'elienne.

[F-G] Fried, D. Goldman, W. Three-dimensional affine
crystallographic groups. Advances in Math. 47 (1983), 1-49.

[F-G-H] Fried, D. Goldman, W. Hirsch, M. Affine manifolds with
nilpotent holonomy. Comment. Math. Helv. 56 (1981) 487-523.

[G1] Goldman, W. Two examples of affine manifolds. Pacific J.
Math. 94 (1981) 327-330.

[G2] Goldman, W. The symplectic nature of fundamental groups of
surfaces. Advances in in Math. 54 (1984) 200-225.

[G3] Goldman, W. Geometric structure on manifolds and varieties of
representations. 169-198, Contemp. Math., 74

[G-H1] Goldman, W. Hirsch, M. The radiance obstruction and
parallel forms on affine manifolds. Trans. Amer. Math. Soc. 286
(1984), 629-949.

[G-H2] Goldman, W. Hirsch, M. Affine manifolds and orbits of
algebraic groups. Trans. Amer. Math. Soc. 295 (1986), 175-198.

[Gr] Grothendieck, A. Pursuing stacks, preprint
 avaible at University of Bangor.

[God] Godement R. Topologie alg\'ebrique et th\'eorie des
faisceaux. (1958) Hermann.

[K1 Koszul, J-L. Vari\'et\'es localement plates et convexit\'e.
Osaka J. Math. (1965), 285-290.

[K2] Koszul, J-L. D\'eformation des connexions localement plates.
Ann. Inst. Fourier 18 (1968), 103-114.

[Mc] Maclane, S. Homology. Springer-Verlag, 1963.

[Ma] Margulis, G. Complete affine locally flat manifolds with a
free fundamental group. J. Soviet. Math. 134 (1987), 129-134.

[Mi] Milnor, J. W. On fundamental groups of complete affinely flat
manifolds, Advances in Math. 25 (1977) 178-187.

[S-T] Sullivan, D. Thurston, W. Manifolds with canonical
coordinate charts: some examples. Enseign. Math 29 (1983), 15-25.

[T1] Tsemo, A. Th\`ese, Universit\'e de Montpellier II. 1999.

[T2] Tsemo, A. Automorphismes polynomiaux des vari\'et\'es
affines. C.R. Acad. Sci. Paris S\'erie I Math 329 (1999) 997-1002.

[T3] Tsemo, A. D\'ecomposition des vari\'et\'es affines. Bull.
Sci. Math. 125 (2001) 71-83.

[T4] Tsemo, A. Dynamique des vari\'et\'es affines. J. London Math.
Soc. 63 (2001) 469-487.

[T5] Tsemo, A. Fibr\'es affines to be published in Michigan J.
Math. vol 49.

[T6] Tsemo, A. Composition series of affine manifolds and
$n-$gerbes Submitted.

[W] Witten, E. Quantum field theory and the Jones Polynomial,
Comment. Math. Phys. 121 (1989) 351-399.

\end{document}